\documentclass{beamer}

\usepackage{amssymb}

\usepackage{beamerthemesplit}

\usepackage{amsmath,amssymb}
\usepackage{enumerate}

\usepackage{amssymb}
\usepackage{amsmath}
\usepackage{color}
\usepackage{graphics}
\usepackage{epsfig}
\usepackage[all,knot]{xy}
\xyoption{arc}

\def\qed{\hfill{$\square$}}

\def\qed{\hfill\vrule width2mm height2mm depth2mm}
\def\v|{\;\vrule width1pt height3mm depth0mm\;}

\def\ex{\begin{example}\rm}
\def\endex{\qed\end{example}}
\def\rem{\begin{rem}\rm}
\def\endrem{\qed\end{rem}}

\def\0{{\bf 0}}
\def\1{{\bf 1}}

\def\cA{{\cal A}}

\def\cF{{\cal F}}

\def\<{\langle}
\def\>{\rangle}
\def\({\left(}
\def\){\right)}

\def\beq#1{\color{blue}\begin{equation}\label{#1}}
\def\eeq{\end{equation}\color{black}}

\title{On piecewise smooth cohomology \\ of locally trivial Lie groupoids}
\author{Jose M. R. Oliveira}
\date{Moscow, May 2012}

\begin{document}

\frame{\titlepage}

\section{Introduction}
\subsection{Piecewise smooth cohomology of transitive Lie algebroids}

\frame{
\frametitle{Piecewise smooth cohomology of transitive Lie algebroids}
\framesubtitle{Mishchenko-Oliveira theorem}
Let $M$ be a smooth manifold, smoothly triangulated by a simplicial complex $K$, and $\cA$ a transitive Lie algebroid on $M$.
\begin{itemize}
\item Since $\cA$ is transitive, the Lie algebroid restriction $\cA_{\Delta}^{!!}$ exists for each simplex $\Delta\in K$.
\item If $\Delta$ and $\Delta'$ are two simplices of $K$ and $\Delta'$ is face of $\Delta$ then $$\big(\cA_{\Delta}^{!!}\big)_{\Delta'}^{!!}=\cA_{\Delta'}^{!!}$$
\end{itemize}
and so the family $(\cA_{\Delta})_{\Delta\in K}$ is a complex of transitive Lie algebroids on $K$.
}

\frame{
\frametitle{Piecewise smooth cohomology of transitive Lie algebroids}
\framesubtitle{Mishchenko-Oliveira theorem}
\begin{itemize}
\item \textbf{Definition - piecewise smooth forms}. Let $M$ be a smooth manifold, smoothly triangulated by a simplicial complex $K$, and $\cA$ a transitive Lie algebroid on $K$. A
piecewise smooth form of degree $p$ ($p\geq 0$) on $\cA$ is a family $(\omega_{\Delta})_{\Delta \in K}$ of smooth forms, each form $\omega_{\Delta}$ defined on $\cA_{\Delta}^{!!}$, which are compatible under restrictions of $\cA$ to faces of simplices of $K$, that is, for each $\Delta$, $\Delta'$ $\in K$, with $\Delta'$ face of
$\Delta$, $$(\omega_{\Delta})^{!!}_{/\Delta'}=\omega_{\Delta'}$$
\item $\Omega^{\ast}(\cA;K)$ denotes the commutative cochain algebra of all piecewise smooth forms on $\cA$, which is defined over $\textrm{\textbf{R}}$. The cohomology space of this complex will be denoted by $H^{\ast}(\cA;K)$.
\end{itemize}
}

\frame{
\frametitle{Piecewise smooth cohomology of transitive Lie algebroids}
\framesubtitle{Mishchenko-Oliveira theorem}
\textbf{Mishchenko-Oliveira theorem}. Let $M$ be a smooth manifold, smoothly triangulated by a simplicial complex $K$, and $\cA$ a transitive Lie algebroid on $K$. The map $$r:\Omega^{\ast}(\cA;M)\longrightarrow \Omega^{\ast}(\cA;K)$$ given by $$r(\omega)=(\omega_{\Delta}^{!!})_{\Delta \in K}$$ induces an isomorphism in cohomology.
}

\frame{
\frametitle{Piecewise smooth cohomology of transitive Lie algebroids}
\framesubtitle{Mishchenko-Oliveira theorem}
\textbf{Some references}.

\begin{itemize}
\item Hassler Whitney - Geometric Integration Theory, chapter VII (Princeton
University Press, 1957).
\item Dennis Sullivan - Infinitesimal computations in topology, chapter 7 (Publ.
I.H.E.S. 47 (1977) 269-331).
\item Kirill Mackenzie - General Theory of Lie Groupoids and Lie Algebroids (London Mathematical Society Lecture Note Series 213, Cambridge U. Press, 2005).
\item Jan Kubarski - The Chern-Weil homomorphism of regular Lie
algebroids (Publ. Dep. Math. University of Lyon 1, 1991).
\end{itemize}
}

\section{Lie groupoids - introduction}
\subsection{Restrictions of Lie groupoids}

\frame{
\frametitle{Restrictions of Lie groupoids}
\framesubtitle{Lie algebroid of a Lie groupoid}
Let $M$ be a smooth manifold and $G$ a Lie groupoid on $M$ with source projection $\alpha:G\longrightarrow M$ and target projection $\beta:G\longrightarrow M$. Denote by $1:M\longrightarrow G$ the object inclusion map of $G$ and $G_{x}=\alpha^{-1}(x)$ the $\alpha$-fibre of G in $x$, for each $x\in M$. The Lie algebroid of $G$ is $(\cA(G),[\cdot,\cdot],\gamma)$, in which

\begin{itemize}
\item $\cA(G)=\bigsqcup _{x\in M}T_{1_{x}}G_{x}$ (disjoint union).
\item The anchor $\gamma:\cA(G)\longrightarrow TM$ is defined by $\gamma(a)=D\beta_{1_{x}}(a)$.
\item The Lie bracket is defined by $[\xi,\eta]=[\xi',\eta']_{G}$, where $\xi'_{o}$ and $\eta'_{o}$ denote the unique $\alpha$-right-invariant vector fields on $G$ such that $\xi'_{1_{x}}=\xi_{x}$ and $\eta'_{1_{x}}=\eta_{x}$, $\forall x\in M$ (see Kubarski's paper).
\end{itemize}
}

\frame{
\frametitle{Restrictions of Lie groupoids}
\framesubtitle{Locally trivial Lie groupoids}
Let $M$ be a smooth manifold and $G$ a Lie groupoid on $M$ with source projection $\alpha:G\longrightarrow M$ and target projection $\beta:G\longrightarrow M$.
\textbf{Definition}.
\begin{itemize}
\item The anchor of $G$ is the map $(\beta,\alpha):G\longrightarrow M\times M$.
\item $G$ is called locally trivial if its anchor $(\beta,\alpha):G\longrightarrow M\times M$ is a surjective submersion.
\end{itemize}
}

\frame{
\frametitle{Restrictions of Lie groupoids}
\framesubtitle{Locally trivial Lie groupoids}
\textbf{Proposition}. Let $G$ be a Lie groupoid on a smooth manifold $M$ with source projection $\alpha:G\longrightarrow M$ and target projection $\beta:G\longrightarrow M$.
\begin{itemize}
\item If $G$ is locally trivial then $\cA(G)$ is transitive Lie algebroid (see Mackenzie's book).
\end{itemize}
}

\frame{
\frametitle{Restrictions of Lie groupoids}
\framesubtitle{Transversal submanifolds}
\textbf{Definition}. Let $M$ be a smooth manifold and $G$ a Lie groupoid on $M$ with source projection $\alpha:G\longrightarrow M$ and target projection $\beta:G\longrightarrow M$. Let $N$ be a submanifold of $M$. 
\begin{itemize}
\item The vertex group of $G$ at $N$ is the set $$G^{N}_{N}:=\alpha^{-1}(N)\cap \beta^{-1}(N)$$
\item The submanifold $N$ is called transversal to $G$ if the anchor $(\beta,\alpha):G\longrightarrow M\times M$ and the map $i\times i: N\times N\longrightarrow M\times M$ are transversal.
\end{itemize}
}

\frame{
\frametitle{Restrictions of Lie groupoids}
\framesubtitle{Transversal submanifolds}
\textbf{Proposition}. Let $M$ be a smooth manifold and $G$ a Lie groupoid on $M$ with source projection $\alpha:G\longrightarrow M$ and target projection $\beta:G\longrightarrow M$.
\begin{itemize}
\item If a submanifold $N$ of $M$ is transversal to the Lie groupoid $G$, then the vertex $G^{N}_{N}$ of $G$ at $N$ is a Lie subgroupoid of $G$ with base $N$ and $G^{N}_{N}$ is called the Lie groupoid restriction of $G$ to $N$.
\item If $G$ is locally trivial then $$\cA(G^{N}_{N})=(\cA(G))^{!!}_{\Delta}$$
\end{itemize}
}

\frame{
\frametitle{Restrictions of Lie groupoids}
\framesubtitle{Transversal submanifolds}
\textbf{Proposition}. Let $M$ be a smooth manifold, smoothly triangulated by simplicial complex $K$, and $G$ a Lie groupoid on $M$ with source projection $\alpha:G\longrightarrow M$ and target projection $\beta:G\longrightarrow M$. 
\begin{itemize}
\item Each simplex $\Delta\in K$ is transversal to $G$. Consequently, $\cA(G^{\Delta}_{\Delta})\simeq \cA(G)_{\Delta}^{!!}$.
\item The family $\{\cA(G^{\Delta}_{\Delta})\}_{\Delta\in K}$ is a complex of Lie algebroids on $M$.
\end{itemize}
}

\section{de Rham cohomology of Lie groupoids}
\subsection{Restriction of Lie groupoids}

\frame{
\frametitle{de Rham cohomology of Lie groupoids}
\framesubtitle{Smooth forms on Lie groupoids}
Let $G$ be a Lie groupoid on $M$ with source projection $\alpha:G\longrightarrow M$ and target projection $\beta:G\longrightarrow M$.
\begin{itemize}
\item $\alpha$ is a surjective submersion and so $\alpha$ induces a foliation $\cF$ on $G$. $T \cF$ denotes the tangent bundle of $\cF$.
\item A smooth $\alpha$-form of degree $p$ on the Lie groupoid $G$ is a smooth section of the exterior vector bundle $\bigwedge ^{p}(T^{\ast}\cF; \textbf{R}_{M})$ in which $\textbf{R}_{M}=M\times \textbf{R}$ denotes the trivial vector bundle on $M$. Then, a smooth form is a family $\omega=(\omega_{g})_{g\in G}$ defined on $G$ such that, for each $g\in G$, one has $$\omega_{g}\in \Lambda^{p}\big(\bigsqcup _{g\in G}T_{g}^{\ast}G_{\alpha(g)};\textbf{R}\big)$$
\end{itemize}
}

\frame{
\frametitle{de Rham cohomology of Lie groupoids}
\framesubtitle{Smooth forms on Lie groupoids}
\begin{itemize}
\item The set $\Omega_{\alpha}^{\ast}(G)$ of all smooth $\alpha$ forms is a commutative cochain algebra defined over $\textrm{\textbf{R}}$.
\item The set $\Omega^{p}_{\alpha,L}(G)$ consisting of all $\alpha$-forms on $G$ which are invariant under the groupoid left translations is a cochain subalgebra of $(\Omega_{\alpha}^{\ast}(G),d^{\ast}_{\alpha})$. Its cohomology is denoted by $H^{\ast}_{\alpha,L}(G)$.
\end{itemize}
}

\frame{
\frametitle{de Rham cohomology of Lie groupoids}
\framesubtitle{Smooth forms on Lie groupoids}
\textbf{Proposition}. There is an isomorphism $$\psi:\Omega^{\ast}_{\alpha,L}(G)\longrightarrow \Omega^{\ast}(\cA(G))$$ of cochain algebras defined by $\psi(\omega)_{x}=\omega_{1_{x}}$. Consequently, $$H^{\ast}_{\alpha,L}(G)\simeq H^{\ast}(\cA(G))$$
}

\section{Piecewise de Rham cohomology of Lie groupoids}
\subsection{Main result}

\frame{
\frametitle{Piecewise de Rham cohomology of Lie groupoids}
\framesubtitle{Formulation of the problem}
Let $M$ be a smooth manifold, smoothly triangulated by a simplicial complex $K$, and $G$ a locally trivial Lie groupoid on $M$ with source projection $\alpha:G\longrightarrow M$ and target projection $\beta:G\longrightarrow M$. For each simplex $\Delta\in K$, the Lie algebroid $\cA(G^{\Delta}_{\Delta})$ is transitive and $\cA(G^{\Delta}_{\Delta})=(\cA(G))^{!!}_{\Delta}$. Therefore, we can consider the complex of Lie algebroids $\{\cA(G^{\Delta}_{\Delta})\}_{\Delta\in K}$. Similarly to piecewise smooth forms on Lie algebroids, we give now the notion of piecewise smooth form on $G$.
}

\frame{
\frametitle{Piecewise de Rham cohomology of Lie groupoids}
\framesubtitle{Formulation of the problem - Continued}
\textbf{Definition}. A piecewise smooth $\alpha$-form of degree $p$ ($p\geq 0$) on $G$ is a family $(\omega_{\Delta})_{\Delta \in K}$ such that the
following conditions are satisfied.
\begin{itemize}
\item For each $\Delta\in K$, $\omega_{\Delta}\in \Omega^{p}_{\alpha,L}(G^{\Delta}_{\Delta})$ $\alpha$-smooth form of degree $p$ on $G^{\Delta}_{\Delta}$.
\item If $\Delta$ and $\Delta'$ are two simplices of $K$ such that $\Delta'\prec \Delta$, one has
$(\omega_{\Delta})_{/\Delta'}=\omega_{\Delta'}$.
\end{itemize}
}

\frame{
\frametitle{Piecewise de Rham cohomology of Lie groupoids}
\framesubtitle{Formulation of the problem - Continued}
The $C^{\infty}(G)$-module $\Omega^{\ast}_{\alpha,L,ps}(G)$ of all piecewise $\alpha$-smooth forms of degree $p$ on $G$ is a commutative differential graded algebra defined over $\textrm{\textbf{R}}$. The cohomology space of this complex will be denoted by $H^{\ast}_{\alpha,L,ps}(G)$.
}

\frame{
\frametitle{Piecewise de Rham cohomology of Lie groupoids}
\framesubtitle{Formulation of the problem - Continued}
Our aim is to relate:
\begin{itemize}
\item The cohomology space $H^{\ast}_{\alpha,L,ps}(G)$ of $G$ to the cohomology space $H^{\ast}_{ps}(\cA(G))$ of its Lie algebroid $\cA(G)$.
\item The cohomology space $H^{\ast}_{\alpha,L}(G)$ to the cohomology space $H^{\ast}_{\alpha,L,ps}(G)$.
\end{itemize}
}

\frame{
\frametitle{Piecewise de Rham cohomology of Lie groupoids}
\framesubtitle{Main result}

\textbf{Main result}. Let $M$ be a smooth manifold, smoothly triangulated by a simplicial complex $K$, and $G$ a locally trivial Lie groupoid on $M$ with source projection $\alpha:G\longrightarrow M$ and target projection $\beta:G\longrightarrow M$. Then,

\begin{itemize}
\item $H^{\ast}_{\alpha,L,ps}(G)$ $\simeq $ $H^{\ast}_{ps}(\cA(G))$.
\item $H^{\ast}_{\alpha,L}(G)$ $\simeq $ $H^{\ast}_{\alpha,L,ps}(G)$.
\end{itemize}
}

\frame{
\frametitle{Piecewise de Rham cohomology of Lie groupoids}
\framesubtitle{Main result - Proof}
The first result holds because we can define a map $$\psi:\Omega^{\ast}_{\alpha,L,ps}(G)\longrightarrow \Omega^{\ast}_{ps}(\cA(G))$$ by the condition of $\psi_{\Delta}:\Omega^{p}_{\alpha,L}(G^{\Delta}_{\Delta})\longrightarrow \Omega^{p}(\cA(G^{\Delta}_{\Delta}))$ to be an isomorphism for each simplex $\Delta\in K$.
}

\frame{
\frametitle{Piecewise de Rham cohomology of Lie groupoids}
\framesubtitle{Main result - Proof}
The second result holds because the following diagram
$$
\xymatrix{
\Omega^{\ast}_{\alpha,L}(G)\ar[d]_{r_{G}}\ar[r]^{iso} & \Omega^{\ast}(\cA;M)\ar[d]^{r_{\cA}} \\
 \Omega^{\ast}_{\alpha,L,ps}(G)\ar[r]^{iso} & \Omega^{\ast}_{ps}(\cA;K)
}
$$
is commutative, where $r_{\cA}$ is the restriction map given at Mishchenko-Oliveira's theorem. By this theorem, the map $r_{\cA}$ is an isomorphism in cohomology and so the proof is done.
}

\frame{
\frametitle{Piecewise de Rham cohomology of Lie groupoids}
\framesubtitle{Main result}
\textbf{Corollary}. Our last proposition says that the piecewise de Rham cohomology of a locally trivial Lie groupoid over a combinatorial manifold doesn't depend on the triangulation of the base. Precisely, let $M$ be a smooth manifold smoothly triangulated by a simplicial complex $K$ and $L$ other simplicial complex which a subdivision of $K$. Let $G$ be a locally trivial Lie groupoid on $M$. Then, the piecewise de Rham cohomology of $G$ obtained by the combinatorial manifold $(M,K)$ is isomorphic to the the piecewise de Rham cohomology of $G$ obtained by the combinatorial manifold $(M,L)$. Thus, this isomorphism is induced by the restriction map.
}

\frame{
\frametitle{Piecewise de Rham cohomology of Lie groupoids}
\framesubtitle{Main result}

}


\textbf{References}
\begin{thebibliography}{99}

\bibitem{}
Aleksandr Mishchenko, Ribeiro Jose, \emph{Generalization of the Sullivan construction for
Transitive Lie Algebroids}, available at arXiv:1102.5698v1 [math.AT] 28 Feb 2011.

\bibitem{}
Jan Kubarski, \emph{Pradines-type groupoids over foliations; cohomology, connections and the Chern-Weil homomorphism}, Preprint n. 2, Institute of Mathematics, Technical University of Lodz, August 1986.

\bibitem{}
Kirill Mackenzie, \emph{General Theory of Lie Groupoids and Lie Algebroids}, London Mathematical Society Lecture Note
Series 213, Cambridge U. Press, 2005.

\end{thebibliography}
\end{document}